\documentclass[a5paper]{book}
\usepackage[english,russian]{babel}
\usepackage{inputenc}
\usepackage{amsmath,latexsym,amsthm}
\usepackage{graphicx}
\usepackage{amssymb}

\begin{document}
\begin{flushright}
\begin{tabular}{l}
{\sf Uzbek  Mathematical}\\
{\sf Journal, 2018, 1, pp.\pageref{Jur1}-\pageref{Jur2}}\\
\end{tabular}
\end{flushright}
UDC 517.946.22; AMS Subject Classifications: Primary 35A08
\begin{center}
\textbf{ Fundamental solutions of generalized bi-axially \\ symmetric multivariable Helmholtz equation}\\
\textbf{Ergashev T.G., Hasanov A. }
\end{center}

{\small
\begin{center}
\begin{tabular}{p{9cm}}

Maqolada ko'p o'zgaruvchili umumlashgan ikki o'qli simmetrik
Gelmgolts tenglamasining to'rtta fundamental yechimlari oshkor
ko'rinishda topilgan va ular yangi kiritilgan uch o'zgaruvchili
konflyuent gipergeometrik funksiya yordamida ifodalangan. Topilgan
fundamental yechimlarning maxsuslik tartibi aniqlangan hamda
elliptik tenglamalar uchun  chegaraviy masalalarni yechishda zarur
bo'ladigan boshqa bir necha xossalari o'rganilgan.
%This will be done by editors
\\[0.5 cm]

 Основным результатом настоящей работы является
построение четырех фундаментальных решений обобщенного
двуосесимметрического многомерного уравнения Гельмгольца в явном
виде, которые удалось выразить через недавно введенную вырожденную
гипергеометрическую функцию от трех переменных. Кроме того,
определен порядок особенности и выявлены  свойства найденных
фундаментальных решений, необходимые при решении краевых задач для
вырождающихся эллиптических уравнений второго порядка.

\end{tabular}\end{center} }
\vspace{0.5 cm}  \makeatletter
\renewcommand{\@evenhead}{\vbox{\thepage \hfil {\it Ergashev T.G., Hasanov A.}   \hrule }}
\renewcommand{\@oddhead}{\vbox{\hfill
{\it Fundamental solutions of generalized bi-axially symmetric ...
}\hfill \thepage \hrule}} \makeatother

\label{Jur1}

\textbf{1.Introduction}

It is known that fundamental solutions have an essential role in
studying partial differential equations. Formulation and solving
of many local and non-local boundary value problems are based on
these solutions. Moreover,\,\,\,fundamental solutions appear as
potentials, for instance, as simple-layer and double-layer
potentials in the theory of potentials.

The explicit form of fundamental solutions gives a possibility to
study the considered equation in detail. For example, in the works
of Barros-Neto and Gelfand [1], fundamental solutions for Tricomi
operator, relative to an arbitrary point in the plane were
explicitly calculated. We also mention Leray's work [2], which it
was described as a general method, based upon the theory of
analytic functions of several complex variables, for finding
fundamental solutions for a class of hyperbolic linear
differential operators with analytic coefficients. Among other
results in this direction, we note a work by Itagaki [3], where 3D
high-order fundamental solutions for a modified Helmholtz equation
were found. The found solutions can be applied with the boundary
particle method to some 2D inhomogeneous problems, for example,
see [4].

Singular partial differential equations appear at studying various
problems of aerodynamics and gas dynamics [5]  and irrigation
problems [6]. For instance, the famous Chaplygin equation [7]
describes subsonic, sonic and supersonic flows of gas. The theory
of singular partial differential equations has many applications
and possibilities of various theoretical generalizations. It is,
in fact, one of the rapidly developing branches of the theory of
partial differential equations.

In most cases boundary value problems for singular partial
differential equations are based on fundamental solutions for
these equations, for instance, see [8].

Let us consider the generalized bi-axially symmetric Helmholtz
equation with $p$ variables
$$ H_{\alpha,\beta}^{p,\lambda}(u)\equiv\sum\limits_{i=1}^p \frac{\partial^2u}{\partial x_i^2}+\frac{2\alpha}{x_1}\frac{\partial u}{\partial x_1}+\frac{2\beta}{x_2}\frac{\partial u}{\partial x_2}-\lambda^2u=0 \eqno(1.1) $$
in the domain $R_p^+\equiv
\left\{(x_1,...,x_p):x_1>0,x_2>0\right\}$ , where $p$ is a
dimension of a Euclidean space $(p\geq2)$, $\alpha,\beta$ and
$\lambda$ are constants and $0<2\alpha,2\beta<1$.

In the article [9], the equation (1.1) was considered in two
cases: (1) when $p=2, \alpha=0,\beta>0;$ (2) when
$p=2,\lambda=0,\beta>0.$ In the work [10] in order to find
fundamental solutions, at first two new confluent hypergeometric
functions  were introduced. Furthermore, by means of the
introduced hypergeometric function fundamental solutions of the
equation(1.1) were constructed in an explicit form. For studying
the properties of the fundamental solutions, the introduced
confluent hypergeometric functions are expanded in products by
Gauss's hypergeometric functions. The logarithmic singularity of
the constructed fundamental solutions of equation (1.1) was
explored with the help of the obtained expansion.Fundamental
solutions of equation (1.1) with $p=3$ \,and $\lambda=0$ were used
in the investigation of the Dirichlet problem for
three-dimensional elliptic equation with two singular coefficients
[11].

In the present article for the equation (1.1) in the domain
$R_p^+$ at $p>2$ four fundamental solutions are constructed in
explicit form. Furthermore, some properties of these solutions are
shown, which will be used for solving boundary value problems for
aforementioned equation.

 \textbf{2.About one confluent hypergeometric function of three variables}

The confluent hypergeometric function of three variables which we
will use in the present work looks like [10]
$$A_2(a;b_1,b_2;c_1,c_2;x,y,z)=\sum\limits_{m,n,k=0}^\infty\frac{(a)_{m+n-k}(b_1)_m(b_2)_n}{(c_1)_m(c_2)_nm!n!k!}x^my^nz^k, \eqno(2.1)$$
where $a,\,b_1,\,b_2,\,c_1,\,c_2$ are complex constants,\,\,
$c_1,\,c_2\neq 0,-1,-2,...$ and \\ $(a)_n=\Gamma(a+n)/\Gamma(a)$
is the Pochhammer symbol.

Using the formula of derivation
$$\frac{\partial^{i+j+k}}{\partial x^i\partial y^j \partial z^k}A_2(a;b_1,b_2;c_1,c_2;x,y,z)=$$
$$=\frac{(a)_{i+j-k}(b_1)_i(b_2)_j}{(c_1)_i(c_2)_j}A_2(a+i+j-k;b_1+i,b_2+j;c_1+i,c_2+j;x,y,z),\eqno (2.2)$$
it is easy to show that the hypergeometric function
$A_2(a;b_1,b_2;c_1,c_2;x,y,z)$ satisfies the system of
hypergeometric equations [10]
 $$\left\{ \begin{matrix}
   x(1-x)\omega_{xx}-xy\omega_{xy}+xz\omega_{xz}+[c_1-(a+b_1+1)x]\omega_x-\\-b_1y\omega_y+b_1z\omega_z-ab_1\omega=0,   \\
   y(1-y)\omega_{yy}-xy\omega_{xy}+yz\omega_{yz}+[c_2-(a+b_2+1)y]\omega_y-\\-b_2x\omega_x+b_2z\omega_z-ab_2\omega=0,   \\
   z\omega_{zz}-x\omega_{xz}-y\omega_{yz}+(1-a)\omega_z+\omega=0,   \\
\end{matrix} \right.\eqno(2.3)$$
where $$\omega(x,y,z)=A_2(a;b_1,b_2;c_1,c_2;x,y,z).\eqno (2.4)$$

Really, by virtue of the derivation formula (2.2), it is easy to
calculate the following expressions
$$\omega_x=\sum\limits_{m,n,k=0}^\infty\frac{(a)_{m+n-k}(b_1)_m(b_2)_n}{(c_1)_m(c_2)_nm!n!k!}\frac{(a+m+n-k)(b_1+m)}{(c_1+m)}x^my^nz^k, \eqno(2.5)$$
$$x\omega_x=\sum\limits_{m,n,k=0}^\infty\frac{(a)_{m+n-k}(b_1)_m(b_2)_n}{(c_1)_m(c_2)_nm!n!k!}\frac{m}{1}x^my^nz^k, \eqno(2.6)$$
$$y\omega_y=\sum\limits_{m,n,k=0}^\infty\frac{(a)_{m+n-k}(b_1)_m(b_2)_n}{(c_1)_m(c_2)_nm!n!k!}\frac{n}{1}x^my^nz^k, \eqno(2.7)$$
$$z\omega_z=\sum\limits_{m,n,k=0}^\infty\frac{(a)_{m+n-k}(b_1)_m(b_2)_n}{(c_1)_m(c_2)_nm!n!k!}\frac{k}{1}x^my^nz^k, \eqno(2.8)$$
$$xy\omega_{xy}=\sum\limits_{m,n,k=0}^\infty\frac{(a)_{m+n-k}(b_1)_m(b_2)_n}{(c_1)_m(c_2)_nm!n!k!}\frac{mn}{1}x^my^nz^k, \eqno(2.9)$$
$$xz\omega_{xz}=\sum\limits_{m,n,k=0}^\infty\frac{(a)_{m+n-k}(b_1)_m(b_2)_n}{(c_1)_m(c_2)_nm!n!k!}\frac{mk}{1}x^my^nz^k, \eqno(2.10)$$
$$x^2\omega_{xx}=\sum\limits_{m,n,k=0}^\infty\frac{(a)_{m+n-k}(b_1)_m(b_2)_n}{(c_1)_m(c_2)_nm!n!k!}\frac{(m-1)m}{1}x^my^nz^k. \eqno(2.11)$$

Substituting (2.5)-(2.11) into the first equation of the system
(2.3), we are convinced that function $\omega(x,y,z)$ satisfies
this equation. We are similarly convinced that function
$\omega(x,y,z)$ satisfies the second and third equations of the
system (2.3).

Having substituted $\omega=x^\tau y^\nu z^\mu \psi(x,y,z)$ in the
system (2.3), it is possible to be convinced that for the values

$\tau:
\,\,\,\,\,\,0,\,\,\,\,\,\,\,\,\,1-c_1,\,\,\,\,\,\,0,\,\,\,\,\,\,\,\,\,1-c_1,$

$\nu:\,\,\,\,\,\,0,\,\,\,\,\,\,\,\,\,\,\,\,\,\,\,0,\,\,\,\,\,\,1-c_2,\,\,\,1-c_2,$

$\mu:\,\,\,\,\,0,\,\,\,\,\,\,\,\,\,\,\,\,\,\,\,0,\,\,\,\,\,\,\,\,\,\,\,\,\,\,0,\,\,\,\,\,\,\,\,\,\,\,\,\,0,$
\\the system has four linearly independent solutions
$$\omega_1(x,y,z)=A_2(a;b_1,b_2;c_1,c_2;x,y,z),\eqno (2.12)$$
$$\omega_2(x,y,z)=x^{1-c_1}A_2(a+1-c_1;b_1+1-c_1,b_2;2-c_1,c_2;x,y,z),\eqno (2.13)$$
$$\omega_3(x,y,z)=y^{1-c_2}A_2(a+1-c_2;b_1,b_2+1-c_2;c_1,2-c_2;x,y,z),\eqno (2.14)$$
\,\,\,\,\,\,\,\,\,$\omega_4(x,y,z)=x^{1-c_1}y^{1-c_2}\times$
$$\times
A_2(a+2-c_1-c_2;b_1+1-c_1,b_2+1-c_2;2-c_1,2-c_2;x,y,z).\eqno
(2.15)$$

In order to further study the decomposition properties of the
products by Gauss's hypergeometric functions, we need to know the
same expansions of the function $A_2(a;b_1,b_2;c_1,c_2;x,y,z).$
For this purpose we shall consider the expression
$$A_2(a;b_1,b_2;c_1,c_2;x,y,z)=$$ $$=\sum\limits_{k=0}^\infty\frac{(-z)^k}{(1-a)_kk!}F_2(a-k;b_1,b_2;c_1,c_2;x,y), \eqno (2.16)$$
where
$$ F_2(a;b_1,b_2;c_1,c_2;x,y)=
\sum\limits_{m,n=0}^\infty\frac{(a)_{m+n}(b_1)_m(b_2)_n}{(c_1)_m(c_2)_nm!n!}x^my^n.$$

In [12] for Appell's hypergeometric function
$F_2(a;b_1,b_2;c_1,c_2;x,y)$ following expansion was found

$$F_2(a;b_1,b_2;c_1,c_2;x,y)=\sum\limits_{i=0}^\infty\frac{(a)_i(b_1)_i(b_2)_i}{(c_1)_i(c_2)_ii!}x^iy^i\times$$
$$\times F(a+i;b_1+i;c_1+i;x)F(a+i;b_2+i;c_2+i;y),
\eqno (2.17)$$ where
$F(a,b;c;z)=\sum\limits_{n=0}^\infty\frac{(a)_n(b)_n}{(c)_nn!}z^n$
is a hypergeometric function of Gauss.

Considering expansion (2.17), from the identity (2.16) we find
[10]

$$A_2(a;b_1,b_2;c_1,c_2;x,y,z)=
\sum\limits_{i,j=0}^\infty\frac{(a)_{i-j}(b_1)_i(b_2)_i}{(c_1)_i(c_2)_ii!j!}x^iy^iz^j
\times$$
$$\times F(a+i-j,b_1+i;c_1+i;x)F(a+i-j,b_2+i;c_2+i;y).\eqno (2.18)$$

By virtue of the formula
$$F(a,b;c;x)=(1-x)^{-b}F\left(c-a,b;c;\frac{x}{x-1}\right),$$
we get from expansion (2.18)
$$A_2(a;b_1,b_2;c_1,c_2;x,y,z)=
(1-x)^{-b_1}(1-y)^{-b_2}\times$$
$$\times\sum\limits_{i,j=0}^\infty\frac{(a)_{i-j}(b_1)_i(b_2)_i}{(c_1)_i(c_2)_ii!j!}\left(\frac{x}{1-x}\right)^i\left(\frac{y}{1-y}\right)^iz^j
\times$$
$$\times F\left(c_1-a+j,b_1+i;c_1+i;\frac{x}{x-1}\right)\times$$ $$\times F\left(c_2-a+j,b_2+i;c_2+i;\frac{y}{y-1}\right).\eqno (2.19)$$
Expansion (2.19) will be used for studying properties of the
fundamental solutions.

We note, that the expansions for the hypergeometric function of
Lauricella  $F_A^{(s)}$ were found in [13].

\textbf{3.Fundamental solutions}

We consider the generalized bi-axially symmetric multivariable
Helmholtz equation in the domain $R_p^+$.\,The equation (1.1) has
the following constructive formulas
$$H^{p,\lambda}_{\alpha,\beta}\left(x_1^{1-2\alpha} u \right)\equiv x_1^{1-2\alpha}H^{p,\lambda}_{1-\alpha,\beta}(u), \eqno(3.1) $$
$$H^{p,\lambda}_{\alpha,\beta}\left(x_2^{1-2\beta} u \right)\equiv x_2^{1-2\beta}H^{p,\lambda}_{\alpha,1-\beta}(u). \eqno(3.2) $$
The constructive formulas (3.1),(3.2) give the possibility to
solve boundary value problems for equation (1.1) for various
values of the parameters $\alpha,\beta$.

We search the solution of equation (1.1) in the form
$$u(x)=P(r)\omega(\xi,\eta,\zeta), \eqno (3.3)$$
where
$$r^2=\sum \limits_{i=1}^p \left(x_i-x_{0i}\right)^2,\,\,\,\,r^2_1=\left(x_1+x_{01}\right)^2+\sum \limits_{i=2}^p\left(x_i-x_{0i}\right)^2,$$
$$r^2_2=\left(x_1-x_{01}\right)^2+\left(x_2+x_{02}\right)^2+\sum \limits_{i=3}^p\left(x_i-x_{0i}\right)^2,$$
$$\xi=\frac{r^2-r^2_1}{r^2}=-\frac{4x_1x_{01}}{r^2},\,\,\,\,
\eta=\frac{r^2-r^2_2}{r^2}=-\frac{4x_2x_{02}}{r^2},$$
$$\zeta=-\frac{\lambda^2}{4}r^2,\,\,\,\, P(r)=\left(r^2\right)^{1-\alpha-\beta-\frac{p}{2}}.$$

Substituting (3.3) into equation (1.1), we have
$$A_1\omega_{\xi\xi}+A_2\omega_{\eta\eta}+A_3\omega_{\zeta\zeta}+B_1\omega_{\xi\eta}+B_2\omega_{\xi\zeta}+B_3\omega_{\eta\zeta}+$$ $$+C_1\omega_\xi+C_2\omega_\eta+C_3\omega_\zeta+D\omega=0,\eqno(3.5)$$
where
$$A_1=P\sum\limits_{i=1}^p \left(\frac{\partial\xi}{\partial x_i}\right)^2,\,\,\,\,\,\,A_2=P\sum\limits_{i=1}^p \left(\frac{\partial\eta}{\partial x_i}\right)^2,\,\,\,\,\,A_3=P\sum\limits_{i=1}^p \left(\frac{\partial\zeta}{\partial x_i}\right)^2,$$
$$B_1=2P\sum\limits_{i=1}^p \frac{\partial\xi}{\partial x_i}\frac{\partial\eta}{\partial x_i},\,\,\,\,B_2=2P\sum\limits_{i=1}^p \frac{\partial\xi}{\partial x_i}\frac{\partial\zeta}{\partial x_i},\,\,\,\,B_3=2P\sum\limits_{i=1}^p \frac{\partial\eta}{\partial x_i}\frac{\partial \zeta}{\partial x_i},$$
$$C_1=2\sum\limits_{i=1}^p \frac{\partial P}{\partial x_i}\frac{\partial \xi}{\partial x_i}+P\sum\limits_{i=1}^p \frac{\partial^2\xi}{\partial x_i^2}+P\left(\frac{2\alpha}{x_1}\frac{\partial \xi}{\partial x_1}+\frac{2\beta}{x_2}\frac{\partial \xi}{\partial x_2} \right),$$
$$C_2=2\sum\limits_{i=1}^p \frac{\partial P}{\partial x_i}\frac{\partial \eta}{\partial x_i}+P\sum\limits_{i=1}^p \frac{\partial^2\eta}{\partial x_i^2}+P\left(\frac{2\alpha}{x_1}\frac{\partial \eta}{\partial x_1}+\frac{2\beta}{x_2}\frac{\partial \eta}{\partial x_2} \right),$$
$$C_3=2\sum\limits_{i=1}^p \frac{\partial P}{\partial x_i}\frac{\partial \zeta}{\partial x_i}+P\sum\limits_{i=1}^p \frac{\partial^2\zeta}{\partial x_i^2}+P\left(\frac{2\alpha}{x_1}\frac{\partial \zeta}{\partial x_1}+\frac{2\beta}{x_2}\frac{\partial \zeta}{\partial x_2} \right),$$
$$D=\sum\limits_{i=1}^p \frac{\partial^2 P}{\partial x_i^2}+\frac{2\alpha}{x_1}\frac{\partial P}{\partial x_1}+\frac{2\beta}{x_2}\frac{\partial P}{\partial x_2} -\lambda^2P.$$

After elementary evaluations we find
$$A_1=-\frac{4P}{r^2}\frac{x_{01}}{x_1}\xi(1-\xi),\,\,\,\,\, A_2=-\frac{4P}{r^2}\frac{x_{02}}{x_2}\eta(1-\eta),\eqno(3.6)$$
$$A_3=-\lambda^2P\zeta, \,\,\,\,\,\,B_1=\frac{4P}{r^2}\frac{x_{01}}{x_1}\xi\eta+\frac{4P}{r^2}\frac{x_{02}}{x_2}\xi\eta,\eqno(3.7)$$
$$ B_2=-\frac{4P}{r^2}\frac{x_{01}}{x_1}\xi\zeta+\lambda^2P\xi,\,\,\,\,\,\, B_3=-\frac{4P}{r^2}\frac{x_{02}}{x_2}\eta\zeta+\lambda^2 P\eta,\eqno(3.8)$$
$$C_1=-\frac{4P}{r^2}\frac{x_{01}}{x_1}\left[2\alpha-\left(2\alpha+\beta+\frac{p}{2}\right)\xi\right]+\frac{4P}{r^2}\frac{x_{02}}{x_2}\beta\xi,\eqno(3.9)$$
$$C_2= \frac{4P}{r^2}\frac{x_{01}}{x_1}\alpha\eta-\frac{4P}{r^2}\frac{x_{02}}{x_2}\left[2\beta-\left(\alpha+2\beta+\frac{p}{2}\right)\eta\right],\eqno(3.10)$$
$$C_3=-\frac{4P}{r^2}\frac{x_{01}}{x_1}\alpha\zeta-\frac{4P}{r^2}\frac{x_{02}}{x_2}\beta\zeta-\lambda^2P\left(\frac{p}{2}-\alpha-\beta\right),\eqno(3.11)$$
$$D=\frac{4P}{r^2}\left[\frac{x_{01}}{x_1}\alpha+\frac{x_{02}}{x_2}\beta \right]\left(\alpha+\beta-1+\frac{p}{2}\right)-\lambda^2P.\eqno(3.12)$$
Substituting equalities (3.6)-(3.12) into equation (3.5), we get
the system of hypergeometric equations

$$\left\{ \begin{matrix}
   \xi(1-\xi)\omega_{\xi\xi}-\xi\eta\omega_{\xi\eta}+\xi\zeta\omega_{\xi\zeta}+\left[2\alpha-\left(2\alpha+\beta+\frac{p}{2}\right)\xi\right]\omega_\xi-\\-\alpha\eta\omega_\eta+\alpha\zeta\omega_\zeta-\alpha\left(\alpha+\beta-1+\frac{p}{2}\right)\omega=0,   \\
   \eta(1-\eta)\omega_{\eta\eta}-\xi\eta\omega_{\xi\eta}+\eta\zeta\omega_{\eta\zeta}+\left[2\beta-\left(\alpha+2\beta+\frac{p}{2}\right)\eta\right]\omega_\eta-\\-\beta\xi\omega_\xi+\beta\zeta\omega_\zeta-\beta\left(\alpha+\beta-1+\frac{p}{2}\right)\omega=0,   \\
   \zeta\omega_{\zeta\zeta}-\xi\omega_{\xi\zeta}-\eta\omega_{\eta\zeta}+\left(2-\alpha-\beta-\frac{p}{2}\right)\omega_\zeta+\omega=0.   \\
\end{matrix} \right.\eqno(3.13)$$

Considering the solutions of the system of hypergeometric
equations (2.12)-(2.15), we define
$$\omega_1(\xi,\eta,\zeta)=A_2\left(\alpha+\beta-1+\frac{p}{2}; \alpha,\beta;2\alpha,2\beta;\xi,\eta,\zeta\right),\eqno (3.14)$$
$$\omega_2(\xi,\eta,\zeta)=\xi^{1-2\alpha}A_2\left(-\alpha+\beta+\frac{p}{2}; 1-\alpha,\beta;2-2\alpha,2\beta;\xi,\eta,\zeta\right),\eqno (3.15)$$
$$\omega_3(\xi,\eta,\zeta)=\eta^{1-2\beta}A_2\left(\alpha-\beta+\frac{p}{2}; \alpha,1-\beta;2\alpha,2-2\beta;\xi,\eta,\zeta\right),\eqno (3.16)$$
\,\,\,\,\,\,\,\,\,$\omega_4(\xi,\eta,\zeta)=\xi^{1-2\alpha}\eta^{1-2\beta}\times$
$$\times
A_2\left(1-\alpha-\beta+\frac{p}{2};
1-\alpha,1-\beta;2-2\alpha,2-2\beta;\xi,\eta,\zeta\right).\eqno
(3.17)$$

Substituting the equalities (3.14)-(3.17) into the expression
(3.3), we get some solutions of the equation (1.1)

$q_1(x,x_0)=k_1\left(r^2\right)^{1-\alpha-\beta-\frac{p}{2}}\times$
$$\times A_2\left(\alpha+\beta-1+\frac{p}{2};
\alpha,\beta;2\alpha,2\beta;\xi,\eta,\zeta\right),\eqno (3.18)$$
\,\,\,\,\,\,\,\,\,$q_2(x,x_0)=k_2\left(r^2\right)^{\alpha-\beta-\frac{p}{2}}x_1^{1-2\alpha}x_{01}^{1-2\alpha}\times$
$$\times A_2\left(-\alpha+\beta+\frac{p}{2};
1-\alpha,\beta;2-2\alpha,2\beta;\xi,\eta,\zeta\right),\eqno
(3.19)$$
\,\,\,\,\,\,\,\,\,$q_3(x,x_0)=k_3\left(r^2\right)^{-\alpha+\beta-\frac{p}{2}}x_2^{1-2\beta}x_{02}^{1-2\beta}\times$
$$\times A_2\left(\alpha-\beta+\frac{p}{2};
\alpha,1-\beta;2\alpha,2-2\beta;\xi,\eta,\zeta\right),\eqno
(3.20)$$
\,\,\,\,\,\,\,\,\,$q_4(x,x_0)=k_4\left(r^2\right)^{-1+\alpha+\beta-\frac{p}{2}}x_1^{1-2\alpha}x_{01}^{1-2\alpha}x_2^{1-2\beta}x_{02}^{1-2\beta}\times$
$$\times
A_2\left(1-\alpha-\beta+\frac{p}{2};
1-\alpha,1-\beta;2-2\alpha,2-2\beta;\xi,\eta,\zeta\right),\eqno
(3.21)$$ where $k_1,...k_4$ are constants which will be determined
at solving boundary value problems for equation (1.1). It is easy
to notice that the considered functions (3.18)-(3.21) possess the
properties

$$ \frac{\partial q_1(x,x_0)}{\partial x_1}\mid_{x_1=0}=0,\,\,\,\,\,\,\,\,\frac{\partial q_1(x,x_0)}{\partial x_2}\mid_{x_2=0}=0, \eqno (3.22)$$
$$ q_2(x,x_0)\mid_{x_1=0}=0,\,\,\,\,\,\,\,\,\,\,\,\frac{\partial q_2(x,x_0)}{\partial x_2}\mid_{x_2=0}=0, \eqno (3.23)$$
$$ \frac{\partial q_3(x,x_0)}{\partial x_1}\mid_{x_1=0}=0,\,\,\,\,\,\,\,\,q_3(x,x_0)\mid_{x_2=0}=0, \eqno (3.24)$$
$$ q_4(x,x_0)\mid_{x_1=0}=0,\,\,\,\,\,\,\,\,q_4(x,x_0)\mid_{x_2=0}=0. \eqno (3.25)$$

From the expansion  (2.19) follows  that the fundamental solutions
(3.18)-(3.21)\,\,at $ r \rightarrow 0 $  possess a singularity of
the order $\frac{1}{r^{p-2}}$, where $p>2$.

\textbf{References}
\begin{enumerate}
\item Barros-Neto J.J., Gelfand I.M. Fundamental solutions for the
Tricomi operator I,II,III, Duke Math.J. 98(3),1999. P.465-483;
111(3),2001.P.561-584; 128(1)\,2005.\,P.119-140.
\item Leray J. Un prolongementa de la transformation de Laplace
qui transforme la solution unitaires d'un opereteur hyperbolique
en sa solution elementaire (probleme de Cauchy,IV),
Bull.Soc.Math.France 90, 1962. P.39-156.
\item Itagaki M. Higher order three-dimensional fundamental
solutions to the Helmholtz and the modified Helmholtz equations.
Eng.\,Anal.\,Bound.\,Elem.\,15,1995. P.289-293.
\item Golberg M.A., Chen C.S. The method of fundamental solutions
for potential, Helmholtz and diffusion problems, in: Golberg
M.A.(Ed.), Boundary Integral Methods-Numerical and Mathematical
Aspects, Comput.Mech.Publ.,1998. P.103-176.
\item Bers L. Mathematical aspects of subsonic and transonic gas
dynamics, New York,London. 1958.
\item Serbina L.I. A problem for the linearized Boussinesq equation with a nonlocal Samarskii condition, Differ.Equ. 38(8),2002. P. 1187-1194.
\item Chapligin S.A. On gas streams, Dissertation, Moscow, 1902 (in
Russian).
\item Salakhitdinov M.S., Hasanov A. A solution of the
Neumann-Dirichlet boundary-value problem for generalized
bi-axially symmetric Helmholtz equation. Complex Variables and
Elliptic Equations. 53 (4), 2008. P.355-364.
\item Marichev O.I. Integral representation of solutions of the
generalized double axial-symmetric Helmholtz equation (in
Russian). Differencial'nye Uravneniya, Minsk,
14(10),1978.P.1824-1831.
\item Hasanov A. Fundamental solutions of the generalized bi-axially symmetric Helmholtz
equation. Complex Variables and Elliptic Equations. Vol.52, No 8,
2007. P. 673-683.
\item Karimov E.T., Nieto J.J. The Dirichlet problem for a 3D
elliptic equation with two singular coefficients. Computers and
Mathematics with Applications. 62, 2011. P.214-224.
\item Burchnall J.L., Chaundy T.W. Expansions of Appell's double
hypergeometric functions. The Quarterly Journal of Mathematics,
Oxford, Ser.12,1941. P.112-128.
\item Hasanov A., Srivastava H.M. Some decomposition formulas
associated with the Lauricella function $F_A^{r}$ and other
multiple hypergeometric functions. Applied Mathematics Letters,
19(2), 2006. P.113-121.

\end{enumerate}

\begin{tabular}{p{7cm}l}
Institute of Mathematics, Tashkent, Uzbekistan \\
\end{tabular}

\textbf{Ko'p o'zgaruvchili umumlashgan ikki o'qli simmetrik \\
Gelmgolts tenglamasining fundamental yechimlari}
\\ \textbf{Ergashev T.G'., Hasanov A.}

\textbf{Фундаментальные решения обобщенного
\\ двуосесимметрического многомерного уравнения Гельмгольца} \\
\textbf{Эргашев Т.Г., Хасанов А.}

\label{Jur2}

% \endinput
\end{document}